\documentclass[amstex,12pt,russian,amssymb]{article}

\usepackage{mathtext}
\usepackage[cp1251]{inputenc}
\usepackage[T2A]{fontenc}
\usepackage[russian]{babel}
\usepackage[dvips]{graphicx}
\usepackage{amsmath}
\usepackage{amssymb}
\usepackage{amsxtra}
\usepackage{latexsym}
\usepackage{ifthen}

\begin{document}

\newcounter{lemma}
\newcommand{\lemma}{\par \refstepcounter{lemma}%
{\bf Лемма \arabic{lemma}.}}

\newcounter{corollary}
\newcommand{\corollary}{\par \refstepcounter{corollary}%
{\bf Следствие \arabic{corollary}.}}

\newcounter{remark}
\newcommand{\remark}{\par \refstepcounter{remark}%
{\bf Замечание \arabic{remark}.}}

\newcounter{theorem}
\newcommand{\theorem}{\par \refstepcounter{theorem}%
{\bf Теорема \arabic{theorem}.}}

\newcounter{proposition}
\newcommand{\proposition}{\par \refstepcounter{proposition}%
{\bf Предложение \arabic{proposition}.}}

\renewcommand{\refname}{\centerline{\bf Список литературы}}

\newcommand{\proof}{{\it Доказательство.\,\,}}

\noindent УДК 517.5

{\bf Е.А.~Севостьянов, С.А.~Скворцов} (Житомирский государственный
университет им.~И.~Франко)

{\bf Є.О.~Севостьянов, С.О.~Скворцов} (Житомирський державний
університет ім.~І.~Франко)

{\bf E.A.~Sevost'yanov, S.A.~Skvortsov} (Zhitomir Ivan Franko State
University)

\medskip
{\bf О нульмерности одного класса отображений в метрических
пространствах}

{\bf Про нульвимірність одного класу відображень в метричних
просторах}

{\bf On lightness of one class of mappings in metric spaces}

\medskip\medskip
Для отображений, удовлетворяющих одной оценке модуля семейств кривых
в метрических пространствах, получен результат о нульмерности
прообраза при отображении. Доказано, что отображения,
удовлетворяющие указанной оценке, нульмерны, как только мажоранта,
отвечающая за искажение семейств кривых, имеет конечное среднее
колебание в каждой точке.

\medskip\medskip
Для відображень, що задовольняють одну оцінку модуля сімей кривих у
метричних просторах, отримано результат  про нульвимірність
прообразу при відображенні. Доведено, що відображення, які
задовольняють вказану оцінку, нульвимірні, як тільки мажоранта, що
відповідає за спотворення сімей кривих, має скінченне середнє
коливання в кожній точці.

\medskip\medskip
For mappings in metric spaces satisfying one inequality with respect
to modulus of families of curves, there is proved a lightness of
preimage under the mapping. It is proved that, the mappings,
satisfying estimate mentioned above, are light, whenever a function
which corresponds to distortion of families of curves under the
mapping, is of finite mean oscillation at every point.

\newpage
{\bf 1. Введение.} Основные определения и обозначения, используемые
в статье, могут быть найдены в монографиях \cite{Re}, \cite{Ri} и
\cite{MRSY}.

Как известно, для отображений с ограниченным искажением по Решетняку
имеют место неравенства вида
\begin{equation}\label{eq1}
M(f(\Gamma)\leqslant K\cdot(M(\Gamma))\,,
\end{equation}
где $M$ -- модуль семейств кривых $\Gamma,$ а $K$ -- коэффициент
квазиконформности (см. \cite{Pol}). Установлено, что указанные
отображения открыты и дискретны, в частности -- нульмерны, см.
\cite[теоремы~6.3 и 6.4, $\S\,~6,$ гл.~II]{Re}. (Отображение
$f:D\rightarrow \overline{{\Bbb R}^n}={\Bbb R}^n\cup\{\infty\}$
называется {\it нульмерным}, если ${\rm dim\,}\{f^{\,-1}(y)\}=0$ для
каждого $y\in \overline{{\Bbb R}^n},$ где ${\rm dim}$ обозначает
топологическую размерность множества, см. \cite{HW}). Если,
напротив, задано неравенство вида (\ref{eq1}), a priori не
являющееся отображением с ограниченным искажением, то по поводу его
дискретности (нульмерности), насколько нам известно, ничего нельзя
утверждать, поскольку до сих пор это является открытой проблемой.
Рассмотрим теперь неравенство, <<противоположное к (\ref{eq1})>>:
пусть у нас вместо (\ref{eq1}) имеет место соотношение
\begin{equation}\label{eq2}
M(\Gamma)\leqslant K\cdot(M(f(\Gamma)))\,.
\end{equation}
Такие неравенства также установлены для отображений с ограниченным
искажением: хорошо известно, что в этом случае $K:=N(f, A)K_O(f),$
где $K_O={\rm vrai\,} \sup K_O(x, f),$ $K_O(x, f)$ -- максимальная
дилатация, а $N(f, A)$ -- максимальная кратность отображения $f$ на
множестве $A,$ где лежит семейство $\Gamma$ (см.
\cite[теорема~3.2]{MRV$_1$}, либо \cite[теорема~6.7, гл.~II]{Ri}).

В случае $n$-мерного евклидового пространства отметим наличие
положительного ответа на вопрос о дискретности отображения,
участвующего (\ref{eq2}). Именно, в сравнительно недавней работе
\cite{Sev$_1$} первым автором настоящей статьи доказан следующий
общий результат: предположим, что при заданной измеримой по Лебегу
функции $Q:D\rightarrow [1, \infty]$ отображение $f:D\rightarrow
{\Bbb R}^n,$ $n\ge 2,$ удовлетворяет оценке вида
$$M(\Gamma )\le \int\limits_{D^{\,\prime}} Q(y)\cdot \rho_*^n
(y)d\mu^{\,\prime}(y)\,,$$ где $\rho_*$ -- произвольная
неотрицательная борелевская функция такая, что
$\int\limits_{\gamma_*}\rho_*(y)ds\ge 1\quad\forall\quad \gamma_*\in
f(\Gamma)\,,$ область $D^{\,\prime}=f(D),$ а $M$ -- конформный
модуль семейства кривых. Тогда отображение $f$ нульмерно (более того
-- открыто и дискретно, см. \cite[Теорема]{Sev$_1$}).

\medskip
Чтобы в каком-то смысле подытожить наши исследования в этом
направлении, мы покажем в настоящей заметке, что упомянутый
результат имеет место не только в ${\Bbb R}^n,$ но также и в более
общих метрических пространствах, регулярных по Альфорсу, в которых
выполнено так называемое $(1, p)$-неравенство Пуанкаре, где $p$ --
некоторое число, которое будет оговорено далее.

Теперь немного об определениях, используемых в тексте. Всюду далее
$\left(X,\,d, \mu\right)$ и $\left(X^{\,\prime},\,d^{\,\prime},
\mu^{\,\prime}\right)$
--- произвольные  метрические пространства с метриками $d$ и $d^{\,\prime},$
наделённые локально конечными борелевскими мерами $\mu$ и
$\mu^{\,\prime},$ и конечными хаусдорфовыми размерностями
$\alpha\geqslant 2$ и $\alpha^{\,\prime}\geqslant 2,$
соответственно. Ниже мы считаем известными определения, связанные с
кривыми в метрическом пространстве, длинами дуг, интегралами,
условиями допустимости и так далее (см. \cite[разд.~13]{MRSY}).
Следуя \cite[раздел 7.22]{He} будем говорить, что борелева функция
$\rho\colon  X\rightarrow [0, \infty]$ является {\it верхним
градиентом} функции $u\colon X\rightarrow {\Bbb R},$ если для всех
спрямляемых кривых $\gamma,$ соединяющих точки $x$ и $y\in X,$
выполняется неравенство $|u(x)-u(y)|\leqslant
\int\limits_{\gamma}\rho\,|dx|,$ где, как обычно,
$\int\limits_{\gamma}\rho\,|dx|$ обозначает линейный интеграл от
функции $\rho$ по кривой $\gamma.$ Будем также говорить, что в
указанном пространстве $X$ выполняется {\it $(1; p)$-неравенство
Пуанкаре,} если найдётся постоянная $C\geqslant 1$ такая, что для
каждого шара $B\subset X,$ произвольной локально ограниченной
непрерывной функции $u\colon X\rightarrow {\Bbb R}$ и любого её
верхнего градиента $\rho$ выполняется следующее неравенство:
$$\frac{1}{\mu(B)}\int\limits_{B}|u-u_B|d\mu(x)\leqslant C\cdot({\rm diam\,}B)\left(\frac{1}{\mu(B)}
\int\limits_{B}\rho^p d\mu(x)\right)^{1/p}\,.$$
Метрическое пространство $(X, d, \mu)$ назовём {\it $n$-регулярным
по Альфорсу,} если при каждом $x_0\in X,$ некоторой постоянной
$C\geqslant 1$ и всех $R<{\rm diam}\,X$
$$\frac{1}{C}R^{n}\leqslant \mu(B(x_0, R))\leqslant CR^{n}\,.$$

\medskip Пусть $G$ -- область в метрическом пространстве $(X,d,\mu)$.
Следуя \cite{IR$_1$}, будем говорить, что функция
$\varphi:G\rightarrow{\Bbb R}$ имеет {\it конечное среднее колебание
в точке $x_{0}\in\overline{G}$}, пишем $\varphi \in FMO(x_{0})$,
если
\begin{equation}\label{eq13.4.111} \overline{\lim\limits_{\varepsilon\rightarrow
0}}\,\, \,\frac{1}{\mu(B(x_{0},\varepsilon))}
\int\limits_{B(x_{0},\varepsilon)}|\varphi(x)-\overline{\varphi}_{\varepsilon}|\,\,d\mu(x)<\infty
\end{equation}
где $$\overline{\varphi}_{\varepsilon}
=\frac{1}{\mu(B(x_{0},\varepsilon))}
\int\limits_{B(x_{0},\varepsilon)}\varphi(x)\,\,d\mu(x)$$ -- среднее
интегральное значение функции $\varphi(x)$ над шаром
$$B(x_0,\varepsilon)=\{x\in G: d(x,x_0)<\varepsilon\}$$ по
отношению к мере $\mu$. Условие (\ref{eq13.4.111}) включает в себя
предположение об интегрируемости функции $\varphi$ по отношению к
мере $\mu$ на множестве $B(x_0,\varepsilon)$ при некотором
$\varepsilon>0$.

Пусть $p\geqslant 1,$ тогда {\it $p$-модулем} семейства кривых
$\Gamma$ в метрическом пространстве $X$ называется величина
$$M_p(\Gamma)=\inf\limits_{\rho \in \,{\rm adm}\,\Gamma}
\int\limits_{X} \rho^p(x)d\mu(x).$$
{\it Областью $D$ в метрическом пространстве $X$} называется
множество $D,$ являющееся линейно связным в $X.$ Пусть $E,$
$F\subset X$ -- произвольные множества. Обозначим через
$\Gamma(E,F,D)$ семейство всех кривых $\gamma:[a,b]\rightarrow X$
которые соединяют $E$ и $F$ в $D\,,$ т.е. $\gamma(a)\in
E,\gamma(b)\in\,F$ и $\gamma(t)\in D$ при $t \in (a, b).$ Пусть $D$
-- область в $X.$ Для $y_0\in f(D)$ и чисел $0<r_1<r_2<\infty$
обозначим
\begin{equation}\label{eq1**}
A(y_0, r_1, r_2)=\left\{ y\,\in\,X^{\,\prime}: r_1<d(y,
y_0)<r_2\right\}\,,\,\,\, S(y_0, r)=\{y\,\in\,X^{\,\prime}: d(y,
y_0)=r\}\,.\end{equation}
Если $f:D\rightarrow X^{\,\prime}$ -- заданное отображение, то для
фиксированного $y_0\in f(D)$ и произвольных $0<r_1<r_2<\infty$
обозначим через $\Gamma(y_0, r_1, r_2)$ семейство всех кривых
$\gamma$ в области $D$ таких, что $f(\gamma)\in \Gamma(S(y_0, r_1),
S(y_0, r_2), A(y_0, r_1, r_2)).$ Предположим, $f(A)$ измеримо в
$X^{\,\prime}$ для любой области $A\subset X,$ тогда при $p, q> 1$
рассмотрим вместо неравенство
\begin{equation} \label{eq2*A}
M_p(\Gamma(y_0, r_1, r_2))\leqslant \int\limits_{f(D)} Q(y)\cdot
\eta^q (y) d\mu^{\,\prime}(y)
\end{equation}
выполненное для любой неотрицательной измеримой по Лебегу функции
$\eta: (r_1,r_2)\rightarrow [0,\infty ]$ такой, что
\begin{equation}\label{eqA2}
\int\limits_{r_1}^{r_2}\eta(r) dr\geqslant 1\,.
\end{equation}

\medskip
Основной результат настоящей статьи заключает в себе следующая

\medskip
\begin{theorem}\label{th3}{\sl\, Предположим, $\left(X,\,d, \mu\right)$ и
$\left(X^{\,\prime},\,d^{\,\prime}, \mu^{\,\prime}\right)$ --
метрические пространства с метриками $d$ и $d^{\,\prime},$
наделённые локально конечными борелевскими мерами $\mu$ и
$\mu^{\,\prime},$ а $G$ и $G^{\,\prime}$ -- области в $X$ и
$X^{\,\prime},$ имеющие конечные хаусдорфовы размерности
$\alpha\geqslant 2$ и $\alpha^{\,\prime}\geqslant 2,$
соответственно. Предположим, кроме того, $G$ -- локально компактное
и локально связное пространство, регулярное по Альфорсу, в котором
выполнено $(1; p)$-неравенство Пуанкаре при некотором $p\in
[\alpha-1, \alpha].$

Пусть $f:G\,\rightarrow\,{G}^{\,\prime}$ -- непрерывное отображение,
сохраняющее измеримость областей. Тогда, если $f$ удовлетворяет
(\ref{eq2*A}) в каждой точке $y_0\in f(G)$ при некотором $q\in (1,
\alpha^{\,\prime}]$ и любой неотрицательной измеримой по Лебегу
функции $\eta: (r_1,r_2)\rightarrow [0,\infty ],$ удовлетворяющей
условию (\ref{eqA2}), и $Q\in FMO$ в каждой точке $y_0\in f(G),$ то
отображение $f$ либо нульмерно, либо постоянно в $G.$ }
\end{theorem}

\medskip
{\bf 2. Формулировка и доказательство основной леммы.} Связный
компакт $C\subset X$ будем называть {\it континуумом}. Говорят, что
семейство кривых $\Gamma_1$ {\it минорируется} семейством
$\Gamma_2\,,$ пишем $\Gamma_1 > \Gamma_2,$
если для каждой кривой $\gamma \in \Gamma_1$ существует подкривая,
которая принадлежит семейству $\Gamma_2.$
В этом случае, $M_p(\Gamma_1) \leqslant M_p(\Gamma_2)$ (см., напр.,
\cite[теорема~1]{Fu}).

\medskip
Справедливо следующее утверждение (см.~\cite[предложение~4.7]{AS}).

\medskip
\begin{proposition}\label{pr2}
{\sl Пусть $X$ --- $\alpha$-регулярное по Альфорсу метрическое
пространство с мерой, в котором выполняется $(1; p)$-неравенство
Пуанкаре, $\alpha\geqslant 1,$ $\alpha-1\leqslant p\leqslant\alpha$
Тогда для произвольных континуумов $E$ и $F,$ содержащихся в шаре
$B(x_0, R),$ и некоторой постоянной $C>0$ выполняется неравенство
$$M_p(\Gamma(E, F, X))\geqslant \frac{1}{C}\cdot\frac{\min\{{\rm diam}\,E, {\rm diam}\,F\}}{R^{1+p-\alpha}}\,.$$ }
\end{proposition}

\medskip
Следующая лемма включает в себя основной результат настоящей работы
в наиболее общей ситуации.

\medskip
\begin{lemma}\label{lem1}
{\sl\,Предположим, $\left(X,\,d, \mu\right)$ и
$\left(X^{\,\prime},\,d^{\,\prime}, \mu^{\,\prime}\right)$ --
метрические пространства с метриками $d$ и $d^{\,\prime},$
наделённые локально конечными борелевскими мерами $\mu$ и
$\mu^{\,\prime},$ а $G$ и $G^{\,\prime}$ -- области в $X$ и
$X^{\,\prime},$ имеющие конечные хаусдорфовы размерности
$\alpha\geqslant 2$ и $\alpha^{\,\prime}\geqslant 2,$
соответственно. Предположим, кроме того, $G$ -- локально компактное
и локально связное пространство, регулярное по Альфорсу, в котором
выполнено $(1; p)$-неравенство Пуанкаре при некотором $p\in
[\alpha-1, \alpha].$

Далее, предположим, что для при некотором $q\in (1, n]$ и каждого
$y_0\in D$ найдётся $\varepsilon_0>0,$ для которого выполнено
соотношение
\begin{equation} \label{eq4!}
\int\limits_{A(y_0, \varepsilon, \varepsilon_0)
}Q(y)\cdot\psi^q(d^{\,\prime}(y, y_0)) \
d\mu^{\,\prime}(y)\,=\,o\left(I^q(\varepsilon, \varepsilon_0)\right)
\end{equation}
для некоторой измеримой по Лебегу функции
$\psi(t):(0,\infty)\rightarrow (0,\infty),$ такой что
\begin{equation} \label{eq5}
0< I(\varepsilon,
\varepsilon_0):=\int\limits_{\varepsilon}^{\varepsilon_0}\psi(t)dt <
\infty
\end{equation}
при всех $\varepsilon\in(0,\varepsilon_0),$ где $A(y_0, \varepsilon,
\varepsilon_0)$ определено в (\ref{eq1**}) при $r_1=\varepsilon,$
$r_2=\varepsilon_0.$

Пусть $f:G\,\rightarrow\,{G}^{\,\prime}$ -- непрерывное отображение,
сохраняющее измеримость областей. Тогда, если $f$ удовлетворяет
(\ref{eq2*A}) в каждой точке $y_0\in f(G)$ при некоторой измеримой
функции $Q:G\rightarrow [0, \infty],$ некотором фиксированном $q\in
(0, \alpha^{\,\prime}]$ и любой неотрицательной измеримой по Лебегу
функции $\eta: (r_1,r_2)\rightarrow [0,\infty ],$ удовлетворяющей
условию (\ref{eqA2}), то $f$ либо нульмерно, либо постоянно.}
\end{lemma}

\medskip
\begin{remark}\label{rem1}
В условиях леммы \ref{lem1}, можно считать, что для произвольного
фиксированного $A,$ такого что $0<A<\varepsilon_0,$ и всех
$\varepsilon\in (0, A),$ выполняется условие вида
$\int\limits_{\varepsilon}^{A}\,\psi(t) dt>0.$ Действительно, из
того, что $Q(x)>0$ п.в. и $\psi>0,$ а также соотношений (\ref{eq4!})
и (\ref{eq5}) следует, что $\int\limits_{\varepsilon}^{A}\,\psi(t)
dt\rightarrow \infty$ при $\varepsilon\rightarrow 0,$ поскольку
величина интеграла слева в (\ref{eq4!}) увеличивается при уменьшении
$\varepsilon.$
\end{remark}

\medskip
{\it Доказательство леммы \ref{lem1}.} Если отображение $f$
постоянно, доказывать нечего. Пусть $f\not\equiv const.$ Предположим
противное, а именно, что отображение $f$ не нульмерно. Тогда
найдётся $y_0\in G,$ такое что множество $\{f^{\,-1}(y_0)\}$ не
является всюду разрывным. Следовательно, по определению, существует
невырожденное связное множество $C\subset \{f^{\,-1}(y_0)\}.$
Поскольку пространство $X$ локально компактно, можно считать, что
$C$ -- континуум.

Поскольку по предположению $f\not\equiv y_0,$ ввиду непрерывности
отображения $f$ найдётся $x_0\in G$ и $\delta_0>0:$
$\overline{B(x_0, \delta_0)}\subset G$ и
\begin{equation}\label{eq7*}
f(x)\ne y_0\qquad 
\forall\quad x\in \overline{B(x_0, \delta_0)}\,.
\end{equation}
Ввиду локальной компактности $G$ можно считать, что
$\overline{B(x_0, \delta_0)}$ -- компакт в $X.$ Кроме того, ввиду
локальной связности $G$ найдётся связная окрестность $U\subset
B(x_0, \delta_0).$ Тогда $\overline{U}$ -- континуум в $G.$

\medskip
По предложению \ref{pr2} при $p\in [\alpha-1, \alpha]$ будем иметь
\begin{equation}\label{eq4*}
M_p\left(\Gamma\left(C, \overline{U}, G\right)\right)>0\,.
\end{equation}
Заметим, что в силу неравенства (\ref{eq7*}) и в виду соотношения
$f(C)=\{y_0\},$ ни одна из кривых семейства
$\Delta=f\left(\Gamma\left(C, \overline{U}, G\right)\right)$ не
вырождается в точку. В то же время, все кривые указанного выше
семейства $\Delta$ имеют одним из своих концов точку $y_0.$ Пусть
$\Gamma_i$ -- семейство кривых $\alpha_i(t):(0,1)\rightarrow
G^{\,\prime}$ таких, что $\alpha_i(1)\in S(y_0, r_i),$
$r_i<\varepsilon_0,$ $r_i$ -- некоторая строго положительная
вещественная последовательность, такая что $r_i\rightarrow 0$ при
$i\rightarrow \infty$ и $\alpha_i(t)\rightarrow y_0$ при
$t\rightarrow 0.$  Тогда мы вправе записать:
\begin{equation}\label{eq12*}
\Gamma\left(C, \overline{U}, G\right) =
\bigcup\limits_{i=1}^\infty\,\, \Gamma_i^*\,,
\end{equation}
где $\Gamma_i^*$ -- подсемейство всех кривых $\gamma$ из
$\Gamma\left(C, \overline{U}, G\right),$ таких что $f(\gamma)$ имеет
подкривую в $\Gamma_i.$ Заметим, что для произвольного
$\varepsilon\in (0, r_i)$ ввиду \cite[предложение~13.3]{MRSY}
\begin{equation}\label{eq8*}
\Gamma_i^*>\Gamma(y_0, \varepsilon, r_i)\,.
\end{equation}
Рассмотрим следующую функцию
$$\eta_{i,\varepsilon}(t)=\left\{
\begin{array}{rr}
\psi(t)/I(\varepsilon, r_i), & t\in (\varepsilon, r_i)\,,\\
0,  &  t\not\in (\varepsilon, r_i)\,,
\end{array}
\right. $$
где $I(\varepsilon, r_i)=\int\limits_{\varepsilon}^{r_i}\,\psi (t)
dt.$ Заметим, что
$\int\limits_{\varepsilon}^{r_i}\eta_{i,\varepsilon}(t)dt=1,$
поэтому мы можем воспользоваться неравенством вида (\ref{eq2*A}).
Исходя из этого неравенства, ввиду соотношения (\ref{eq8*}),
получаем:
\begin{equation}\label{eq11*}
M_p(\Gamma_i^*)\leqslant M_p(\Gamma(y_0, \varepsilon, r_i))\leqslant
\int\limits_{A(y_0, \varepsilon, \varepsilon_0)} Q(y)\cdot
\eta^q_{i,\varepsilon}(d^{\,\prime}(y,
y_0))d\mu^{\,\prime}(y)\,\leqslant {\frak F}_i(\varepsilon),
\end{equation}
где
${\frak F}_i(\varepsilon)=\,\frac{1}{{I(\varepsilon,
r_i)}^q}\int\limits_{A(y_0, \varepsilon,
\varepsilon_0)}\,Q(y)\,\psi^{q}(d^{\,\prime}(y,
y_0))\,d\mu^{\,\prime}(y)$ и $I(\varepsilon,
r_i)=\int\limits_{\varepsilon}^{r_i}\,\psi (t) dt.$ Учитывая
$(\ref{eq4!}),$ имеем следующее соотношение:
$$\int\limits_{A(y_0, \varepsilon, \varepsilon_0)}\,Q(y)\,\psi^{q}(d^{\,\prime}(y,
y_0))\,d\mu^{\,\prime}(y)\,=\,
G(\varepsilon)\cdot\left(\int\limits_{\varepsilon}^{\varepsilon_0}\,\psi
(t) dt\right)^q\,,$$
где $G(\varepsilon)\rightarrow 0$ при $\varepsilon\rightarrow 0$ по
условию леммы. Заметим, что
${\frak F}_i(\varepsilon)\,=\,G(\varepsilon)\cdot\left(1\,+\,\frac{
\int\limits_{r_i}^{\varepsilon_0}\,\psi(t)
dt}{\int\limits_{\varepsilon}^{r_i}\,\psi(t) dt}\right)^q,$
где $\int\limits_{r_i}^{\varepsilon_0}\,\psi(t) dt<\infty$ --
фиксированное число, а $\int\limits_{\varepsilon}^{r_i}\,\psi(t)
dt\rightarrow \infty$ при $\varepsilon\rightarrow 0,$ поскольку
величина интеграла слева в $(\ref{eq4!})$ увеличивается при
уменьшении $\varepsilon.$ Таким образом, ${\frak
F}_i(\varepsilon)\rightarrow 0.$ Переходя к пределу в неравенстве
(\ref{eq11*}) при $\varepsilon\rightarrow 0,$ левая часть которого
не зависит от $\varepsilon,$ получаем, что $M_p(\Gamma_i^*)=0$ при
любом натуральном $i.$ Однако, тогда $M_p\left(\Gamma\left(C,
\overline{U}, G\right)\right) =0$ в виду соотношения (\ref{eq12*}) и
того, что
$M_p\left(\bigcup\limits_{i=1}^{\infty}\Gamma_i\right)\leqslant
\sum\limits_{i=1}^{\infty}M_p(\Gamma_i)$ (см. \cite[теорема~1]{Fu}),
что противоречит неравенству (\ref{eq4*}). Полученное противоречие
доказывает, что отображение $f$ является нульмерным, что и
требовалось доказать. $\Box$

\medskip
{\bf 3. О доказательстве основного результата.} Доказательство
теоремы \ref{th3} вытекает из леммы \ref{lem1} и
\cite[лемма~13.2]{MRSY}.

\medskip
Действительно, согласно \cite[лемма~13.2]{MRSY}, условие $Q\in
FMO(y_0)$ влечёт, что при некотором $\varepsilon_0>0$ и
$\varepsilon\rightarrow 0$
$$\int\limits_{\varepsilon<d^{\,\prime}(y, y_0)<\varepsilon_0}
Q(x)\cdot\psi^{q}(d^{\,\prime}(y, y_0))\,
d\mu^{\,\prime}(y)=O\left(\log\log\frac{1}{\varepsilon}\right)\,,$$
где
$\psi(t)=\left(t\,\log{\frac1t}\right)^{-\alpha/q}>0$ и
$0<q\leqslant\alpha.$
Как и прежде, определим $I(\varepsilon,
\varepsilon_0):=\int\limits_{\varepsilon}^{\varepsilon_0}\psi(t)dt,$
тогда
$$
I(\varepsilon,
\varepsilon_0)=\int\limits_{\varepsilon}^{\varepsilon_0}\psi(t) dt
>\log{\frac{\log{\frac{1}
{\varepsilon}}}{\log{\frac{1}{\varepsilon_0}}}}.
$$
В таком случае, заключаем, что условия (\ref{eq4!})--(\ref{eq5})
леммы \ref{lem1} выполнены и, значит, из этой леммы вытекает
требуемое утверждение.~$\Box$

\medskip
В частности, справедливо следующее заключение, вытекающее из леммы
\ref{lem1} и \cite[следствие~13.3]{MRSY}.

\medskip
\begin{corollary}\label{cor1}
{\sl\, Предположим, $\left(X,\,d, \mu\right)$ и
$\left(X^{\,\prime},\,d^{\,\prime}, \mu^{\,\prime}\right)$ --
метрические пространства с метриками $d$ и $d^{\,\prime},$
наделённые локально конечными борелевскими мерами $\mu$ и
$\mu^{\,\prime},$ а $G$ и $G^{\,\prime}$ -- области в $X$ и
$X^{\,\prime},$ имеющие конечные хаусдорфовы размерности
$\alpha\geqslant 2$ и $\alpha^{\,\prime}\geqslant 2,$
соответственно. Предположим, кроме того, $G$ -- локально компактное
и локально связное пространство, регулярное по Альфорсу, в котором
выполнено $(1; p)$-неравенство Пуанкаре при некотором $p\in
[\alpha-1, \alpha].$

Пусть $f:G\,\rightarrow\,{X}^{\,\prime}$ -- непрерывное отображение,
сохраняющее измеримость областей. Тогда, если $f$ удовлетворяет
(\ref{eq2*A}) в каждой точке $y_0\in f(G)$ при некотором $q\in (0,
\alpha^{\,\prime}]$ и любой неотрицательной измеримой по Лебегу
функции $\eta: (r_1,r_2)\rightarrow [0,\infty ],$ удовлетворяющей
условию (\ref{eqA2}), и
$$\limsup\limits_{\varepsilon\rightarrow 0}\frac{1}{\mu^{\,\prime}(B(y_0, \varepsilon))}
\int\limits_{B(y_0, \varepsilon)}Q(y)d\mu^{\,\prime}(y)<\infty$$
 в каждой точке $y_0\in f(G),$ то
отображение $f$ либо нульмерно, либо постоянно в $G.$ }

\end{corollary}

КОНТАКТНАЯ ИНФОРМАЦИЯ

\medskip
\noindent{{\bf Евгений Александрович Севостьянов} \\
\noindent{{\bf Сергей Александрович Скворцов} \\
Житомирский государственный университет им.\ И.~Франко\\
кафедра математического анализа, ул. Большая Бердичевская, 40 \\
г.~Житомир, Украина, 10 008 \\ тел. +38 066 959 50 34 (моб.),
e-mail: esevostyanov2009@mail.ru}

\end{document}